\documentclass{elsart}

\usepackage{amssymb}

\usepackage[english,francais]{babel}

\newtheorem{e-proposition}[theorem]{Proposition}

\newtheorem{e-definition}[theorem]{Definition\rm}

\newtheorem{theoreme}{Th\'eor\`eme}[section]

\newtheorem{remarque}{\it Remarque}

\renewcommand{\theequation}{\arabic{equation}}
\def\og{\leavevmode\raise.3ex\hbox{$\scriptscriptstyle\langle\!\langle$~}}
\def\fg{\leavevmode\raise.3ex\hbox{~$\!\scriptscriptstyle\,\rangle\!\rangle$}}

\def\<{\langle\,}
\def\>{\,\rangle}
\def\cf{{\it cf.$\ $}}

\def\SG{{\mathfrak S}}

\def\Sym{{\bf Sym}}
\def\FQSym{{\bf FQSym}}

\def\MQSym{{\bf MQSym}}

\def\SS{{\bf S}}

\def\N{{\mathbb N}}

\def\sl{\hbox{\goth sl\hskip 1pt}}

\def\K{{\mathbb K}}

\def\MM{{\mathcal M}} 
\def\MMM{{\bf M}}     
           
\def\GG{{\mathcal G}}
\def\supp{{\rm supp}}

\def\GQSym{\hbox{\it GQSym}} 
\def\GTSym{\hbox{\it GTSym}} 
\def\GSym{{\bf GSym}}

\begin{document}

\begin{frontmatter}




%
\selectlanguage{francais}
\title{Alg\`ebres de Hopf de graphes}

\vspace{-2.8cm}
\selectlanguage{english}
\title{Hopf Algebras of Graphs}



\author[mlv]{Jean-Christophe Novelli}
\ead{novelli@univ-mlv.fr}
\author[mlv]{Jean-Yves Thibon}
\ead{jyt@univ-mlv.fr}
\author[mlv,lyon]{Nicolas M. Thi\'ery}
\ead{nthiery@users.sf.net}
\address[mlv]{Institut Gaspard Monge, Universit\'e de
Marne-la-Vall\'ee, 77454 Marne-la-Vall\'ee Cedex 2, France}
\address[lyon]{Laboratoire de Probabilit\'es, Combinatoire et Statistiques,
Universit\'e Lyon I, B\^atiment B, 50, avenue Tony-Garnier, Domaine de
Gerland,  69366 Lyon Cedex 07, France}

\begin{abstract}
We define graded Hopf algebras with bases labeled by various types of
graphs and hypergraphs, provided with natural embeddings into an algebra of polynomials in
infinitely many variables. These algebras are graded by the number of edges
and can be considered as generalizations of symmetric or quasi-symmetric
functions.

\vskip 0.5\baselineskip

\selectlanguage{francais}
\noindent{\bf R\'esum\'e}
\vskip 0.5\baselineskip
\noindent

Nous d\'efinissons des alg\`ebres de Hopf dont les bases sont \'etiquet\'ees
par divers types de graphes et hypergraphes et les r\'ealisons comme sous-alg\`ebres d'une
alg\`ebre de polyn\^omes en une infinit\'e de variables. Ces alg\`ebres sont
gradu\'ees par le nombre d'ar\^etes et peuvent \^etre consid\'er\'ees comme
des g\'en\'eralisations des fonctions sym\'etriques ou quasi-sym\'etriques.
\end{abstract}
\end{frontmatter}


\selectlanguage{francais}
\section{Introduction}
\label{Intro}
On conna\^\i t de nombreux exemples d'alg\`ebres de Hopf gradu\'ees dont les
bases sont index\'ees par des objets combinatoires~\cite{JR,Sch,LRa,NT}.
La notion d'\emph{alg\`ebre de Hopf combinatoire} est utilis\'ee de mani\`ere
informelle depuis quelques ann\'ees, et certains auteurs ont propos\'e de lui
donner un sens pr\'ecis~\cite{ABS}. Si l'on adopte leur d\'efinition, ce sont
les alg\`ebres munies d'un homomorphisme vers l'alg\`ebre des fonctions
quasi-sym\'etriques.
Dans la plupart des exemples connus, l'existence de cet homomorphisme, qui
peut para\^\i tre myst\'erieuse si l'on s'en tient \`a la d\'efinition
abstraite de l'alg\`ebre (dont le produit et le coproduit ne sont souvent
d\'efinis que de mani\`ere r\'ecursive), s'explique simplement par une
r\'ealisation de l'alg\`ebre en termes de polyn\^omes (en variables
commutatives ou non).
Nous nous proposons ici de construire diverses alg\`ebres de Hopf de graphes
directement munies de telles r\'ealisations.

Nous consid\'erons des graphes \'etiquet\'es ou non, orient\'es ou non, avec
ou sans boucles ou ar\^etes multiples, ce qui donne d\'ej\`a 16 alg\`ebres de
Hopf commutatives, mais non cocommutatives.
Par dualit\'e, nous obtenons ensuite 16 alg\`ebres non commutatives.

Il existe d\'ej\`a de nombreux exemples d'alg\`ebres de Hopf
construites sur divers types de graphes (alg\`ebres
d'incidence~\cite{Sch} ou de renormalisation~\cite{Kre,LR1}), mais
celles consid\'er\'ees ici sont d'un type diff\'erent.
%

Nous employons les conventions suivantes pour d\'esigner ces alg\`ebres. Les
acronymes des alg\`ebres commutatives s'\'ecrivent en italiques. Par exemple,
$Sym$ et $QSym$ d\'esignent respectivement les fonctions sym\'etriques et les
fonctions quasi-sym\'etriques.
Les acronymes des alg\`ebres non commutatives sont en caract\`eres gras.
Ainsi, $\Sym$, $\FQSym$, $\MQSym$ correspondent respectivement aux fonctions
sym\'etriques non commutatives, aux fonctions quasi-sym\'etriques libres, et
aux fonctions quasi-sym\'etriques matricielles.
Les noms des alg\`ebres sont choisis en fonction du proc\'ed\'e de
construction. Dans ce qui suit, les alg\`ebres commutatives apparaissent comme
des g\'en\'eralisations naturelles des fonctions quasi-sym\'etriques (graphes
\'etiquet\'es) ou des fonctions sym\'etriques (graphes non \'etiquet\'es).
Elles seront respectivement not\'ees $\GQSym^{\bf v}$ et $\GTSym^{\bf v}$,
${\bf v}$ \'etant un vecteur bool\'een indiquant les options retenues
(orientation, boucles, ar\^etes multiples).
Nous renvoyons le lecteur \`a \cite{NCSF6} pour les autres notations.

Ce travail a b\'en\'efici\'e du soutien du r\'eseau europ\'een ACE
(HPRN-CT-2001-00272).

\section{Graphes \'etiquet\'es, orient\'es, avec boucles et ar\^etes multiples}

Les bases lin\'eaires de l'alg\`ebre $\GQSym^{111}$ sont index\'ees par les
graphes \'etiquet\'es, orient\'es, avec boucles et ar\^etes multiples,
sans sommet isol\'e (c'est-\`a-dire, n'appartenant \`a
aucune ar\^ete), les sommets \'etant num\'erot\'es par des entiers
successifs $1,\dots,m$.
Un tel graphe $G$ \`a $m$ sommets est donc d\'ecrit par une matrice
d'adjacence $A=(a_{ij})_{i,j=1}^m$ telle qu'il n'existe pas d'indice $i$
v\'erifiant $a_{ij}=a_{ji}=0$ pour tout $j$.

Soient $x_{ij}$, $i,j\ge 1$ des ind\'etermin\'ees. Au graphe $G$, nous
associons la s\'erie formelle
\begin{equation}
\MM_G := \sum_{i_1<\ldots<i_m}\prod_{j,k=1}^mx_{i_j i_k}^{a_{jk}}\,.
\end{equation}

Soit $\GQSym^{111}$ le sous-espace vectoriel de $\K[{x_{ij}|i,j\ge 1}]$
engendr\'e par les $\MM_G$, o\`u $\K$ est un corps de caract\'eristique
z\'ero.

Un entier $i\in [0,m]$ est une \emph{coupure admissible} de $G$ s'il n'y a
aucune ar\^ete entre les sommets de $[1,i]$ et ceux de $[i+1,m]$. Nous
notons $C_G$ l'ensemble des coupures admissibles de $G$.

Pour un sous-ensemble de sommets $D\subseteq [1,m]$, nous d\'esignons par
$G|_D$ la restriction du graphe \`a $D$, \emph{avec les sommets
renum\'erot\'es de $1$ \`a $d:=|D|$} en conservant leur ordre initial.
Posons
\begin{equation}
\Delta\MM_G := \sum_{i\in C_G}\MM_{G|_{[1,i]}}\otimes \MM_{G|_{[i+1,m]}}\,.
\end{equation}

\smallskip
\begin{theoreme} {\rm (i)} $\GQSym^{111}$ est une sous-alg\`ebre de
$\K[{x_{ij}|i,j\ge 1}]$. Plus pr\'ecis\'ement, il existe des entiers
$c_{G',G''}^G\in \N$ tels que
\begin{equation}
\MM_{G'}\MM_{G''}=\sum_G c_{G',G''}^G \MM_G\,.
\end{equation}
{\rm (ii)} Le coproduit $\Delta$ est coassociatif et est un morphisme
d'alg\`ebres.
\end{theoreme}
\smallskip

Ainsi, $\GQSym^{111}$ est une alg\`ebre de Hopf gradu\'ee, le degr\'e \'etant
le nombre total d'ar\^etes. Les premi\`eres valeurs des dimensions
$d_n=\dim\GQSym^{111}_n$ sont $1,1,3,39,819,23949$.
On peut les exprimer par une s\'erie infinie
\begin{equation}
\label{dims111}
d_n=\sum_{m\ge 0}\frac{1}{2^{m+1}}{m^2+n-1 \choose n}
\end{equation}
ou comme un produit scalaire de fonctions sym\'etriques,
$d_n=\< h_n\circ h_{11},H_{2n}\>$,
o\`u $H_n$ est le terme de degr\'e $n$ dans la s\'erie
$(1-\sum_{k\ge 1}h_k)^{-1}$ et o\`u $f\circ g$ d\'esigne le pl\'ethysme de $g$
par $f$ (\cf~\cite{McD}).
Il existe des formules analogues pour les autres alg\`ebres de graphes
\'etiquet\'es.

Il est clair que $\Delta$ n'est pas cocommutatif.
De fait, le dual (gradu\'e) $\left(\GQSym^{111}\right)^*$, not\'e
$\GSym^{111}$, est une alg\`ebre libre. Soit $\SS^G=\MM_G^*$ la base duale de
$\MM_G$. Nous dirons qu'un graphe $G$ est \emph{irr\'eductible} s'il n'a pas
de coupure admissible non triviale.

\smallskip
\begin{theoreme}
L'alg\`ebre duale $\GSym^{111}$ est l'alg\`ebre libre sur l'ensemble
$\{\SS^G| G \mbox{\rm\ irr\'eductible}\}$.
\end{theoreme}
\smallskip

Il est possible de munir $\GQSym^{111}$ de plusieurs structures
d'alg\`ebre de Hopf combinatoire, au sens de~\cite{ABS}. Cela revient \`a se
donner un morphisme de Hopf de $\Sym$ vers le dual $\GSym^{111}$.
Soit $\gamma(n)$ le graphe form\'e de $n$ boucles sur un seul sommet (de
matrice d'adjacence $(n)_{1\times 1}$). Il est imm\'ediat que
$S_n\mapsto \SS^{\gamma(n)}$ est un morphisme d'alg\`ebres de Hopf.
On obtient d'autres morphismes de Hopf en rempla\c cant $\gamma(n)$ par les
graphes de matrices \scriptsize$\pmatrix{0&n\cr 0&0}$\normalsize\ ou
\scriptsize$\pmatrix{0&0\cr n&0}$\normalsize.
Si l'on note $\gamma(p,q)$ le graphe de matrice
\scriptsize$\pmatrix{0&p\cr q&0}$\normalsize, l'application
$S_n\mapsto \sum_{p+q=n}\SS^{\gamma(p,q)}$ est encore un morphisme de Hopf.

On remarquera que la sp\'ecialisation $x_{ij}=x_i x_j$ d\'efinit un morphisme
d'alg\`ebres de $\GQSym^{111}$ vers $QSym$, chaque $\MM_G$ \'etant envoy\'ee
sur une $M_I$, mais ce n'est pas un morphisme de cog\`ebres.


\section{Sous-big\`ebres de $\GQSym^{111}$}

Soit $\GQSym^{011}$ le sous-module de $\GQSym^{111}$ engendr\'e par les
$\MM_G$ telles que la matrice d'adjacence de $G$ soit sym\'etrique.
Un tel graphe s'identifie naturellement \`a un graphe non orient\'e.
On d\'efinit de m\^eme $\GQSym^{101}$ comme le sous-module engendr\'e par les
$\MM_G$ telles que la matrice d'adjacence de $G$ ait une diagonale nulle
(graphes sans boucles), et $\GQSym^{001}$ comme l'intersection de
$\GQSym^{011}$ et de $\GQSym^{101}$.

\smallskip
\begin{theoreme}
$\GQSym^{011}$, $\GQSym^{101}$ et $\GQSym^{001}$ sont des sous-alg\`ebres de
Hopf de $\GQSym^{111}$.
\end{theoreme}
\smallskip

On peut aussi d\'efinir des alg\`ebres $\GQSym^{ab0}$ bas\'ees sur des
graphes sans ar\^etes multiples. On ne les obtient pas comme sous-alg\`ebres,
mais comme quotients de $\GQSym^{111}$.

\section{Quotients de $\GQSym^{111}$}

Soit $\GG$ une collection de graphes irr\'eductibles. Soit $\GG^{\sup}$
l'ensemble des graphes dont une restriction $G|_D$ contient toutes les
ar\^etes d'au moins un \'el\'ement de $\GG$.
Alors, les $\MM_G$ pour $G\in\GG^{\sup}$ forment une base d'un id\'eal
(gradu\'e) $I(\GG)$ de $\GQSym^{111}$. C'est \'egalement un co\"\i d\'eal,
de sorte que le quotient $\GQSym^{111}/I(\GG)$ est \`a son tour une alg\`ebre
de Hopf gradu\'ee.
Ce quotient a pour base les classes des $\MM_G$, pour $G\not\in\GG^{\sup}$.

On obtient les graphes sans ar\^etes multiples en prenant pour $\GG$
l'ensemble des trois graphes $\{1\Rightarrow 2\}$, $\{2\Rightarrow 1\}$ et
$\{1\Rightarrow 1\}$, respectivement form\'es d'une ar\^ete double de 1 vers
2, d'une ar\^ete double de 2 vers 1 et d'une boucle double de 1 vers
lui-m\^eme.
Il est facile de retrouver par ce proc\'ed\'e les alg\`ebres
associ\'ees aux graphes sans boucles. On peut aussi obtenir des
alg\`ebres associ\'ees aux partitions non crois\'ees, aux for\^ets,
aux multi-for\^ets et \`a de nombreux autres exemples.


\section{Graphes non \'etiquet\'es}

Les constructions pr\'ec\'edentes s'adaptent au cas des graphes non
\'etiquet\'es. Appelons \emph{support} d'un graphe \'etiquet\'e $G$, le graphe
non \'etiquet\'e sous-jacent $\Gamma$, et notons $\Gamma=\supp(G)$.

\smallskip
\begin{theoreme}
Les sommes
\begin{equation}
\MMM_\Gamma := \sum_{\supp(G)=\Gamma} \MM_G,
\end{equation}
o\`u $\Gamma$ parcourt l'une des classes de graphes non \'etiquet\'es
orient\'es ou non, avec/sans boucles, avec/sans ar\^etes multiples, forment
une base d'une sous-alg\`ebre de Hopf $\GTSym^{abc}$ de $\GQSym^{abc}$.

\end{theoreme}
\smallskip

Les dimensions des composantes homog\`enes de $\GTSym^{abc}$ peuvent se
calculer au moyen des caract\`eres du groupe sym\'etrique. Par exemple, la
dimension de $\GTSym^{111}_n$ est \'egale \`a la multiplicit\'e de la
repr\'esentation triviale de $\SG_m$ dans la $n^{\rm e}$ puissance
sym\'etrique de la repr\'esentation de caract\'eristique
$h_{(m-2,1,1)}+h_{(m-1,1)}$ pour tout $m\geq 2n$.
On obtient la suite $A052171$ de~\cite{Sloane}.

Les alg\`ebres de graphes sur un nombre fini $n$ de sommets
\'etudi\'ees pr\'ec\'edemment par le troisi\`eme auteur~\cite{PT,Thi}
(cf. aussi \cite{Koc})
sont des quotients naturels de $\GTSym^{abc}$, mais ne sont pas de
Hopf.


\smallskip
\begin{remarque}
{\rm
Les constructions pr\'ec\'edentes s'\'etendent \emph{mutatis-mutandis} \`a
des alg\`ebres 
associ\'ees aux hypergraphes $k$-homog\`enes, en prenant cette fois des
variables $x_{i_1,\dots,i_k}$ index\'ees par des $k$-uplets d'entiers.
}
\end{remarque}



%
%
%


\end{document}